\documentclass{amsart}

\usepackage[utf8]{inputenc}
\usepackage[T1]{fontenc}
\usepackage{verbatim}
\usepackage{graphicx}
\usepackage{graphicx,caption2,psfrag,float,color}
\usepackage{amssymb}
\usepackage{amscd}
\usepackage{amsmath}

\usepackage{tikz}
\usetikzlibrary{positioning}

\usepackage{mathabx}

\usepackage[T1]{fontenc}

\newtheorem{theorem}{Theorem}[section]

\newtheorem{lemma}[theorem]{Lemma}
\newtheorem{proposition}[theorem]{Proposition}
\numberwithin{equation}{subsection}
\newtheorem{definition}[theorem]{Definition}

\title{The ternary Goldbach problem with a prime with a missing digit and primes of special types}
\author{Helmut Maier and Michael Th. Rassias}
\date{\today}
\address{Department of Mathematics, University of Ulm, Helmholtzstrasse 18, 89081 Ulm, Germany.}
\email{helmut.maier@uni-ulm.de}
\address{Department of Mathematics and Engineering Sciences,
Hellenic Military Academy,
16673 Vari Attikis, Greece \\
\&
Moscow Institute of Physics and Technology,
141700 Dolgoprudny,
Institutskiy per, d. 9,
Russia \\
\& Institute for Advanced Study, Program in Interdisciplinary Studies,
1 Einstein Dr, Princeton, NJ 08540, USA.}\email{michail.rassias@math.uzh.ch}\thanks{}
\begin{document}

 \maketitle
 
\begin{abstract} 
Let 
$$\gamma^*:=\frac{8}{9}+\frac{2}{3}\:\frac{\log(10/9)}{\log 10}\:(\approx 0.919\ldots)\:,\ \gamma^*<\frac{1}{c_0}\leq 1\:.$$

Let $\gamma^*<\gamma_0\leq 1$, $c_0=1/\gamma_0$ be fixed. Let also $a_0\in\{0,1,\ldots, 9\}$.\\
In \cite{maier-rassias-gold} we proved on assumption of the Generalized Riemann Hypothesis (GRH), that each sufficiently large odd integer $N_0$ can be represented in the form 
$$N_0=p_1+p_2+p_3\:,$$
where for $i=2, 3$ the primes $p_i$ are Piatetski-Shapiro primes -- primes of the form $p_i=[n_i^{c_0}]$, $n_i\in\mathbb{N}$ --  whereas the decimal expansion of $p_1$ does not contain the digit $a_0$.\\
In this paper we replace one of the Piatetski-Shapiro primes $p_2$ and $p_3$ by primes of the type
$$p=x^2+y^2+1\:.$$
\newline
\end{abstract}

\section{Introduction and statement of result}

One of the most famous unsolved problems in number theory is the binary Goldbach problem:\\
\textit{Every even number $\geq 4$ is the sum of two prime numbers.}\\
The ternary Goldbach problem, the representation of odd integers as the sum of three primes, has been treated more successfully (cf. \cite{rasvino}). Let
\[
R(N_0):=\sum_{\substack{(p_1,p_2,p_3)\\ p_1+p_2+p_3=N_0}} (\log p_1)(\log p_2)(\log p_3) \tag{1.1}
\]
(Here and in the sequel the letter $p$ denotes primes).\\
Vinogradov \cite{vino} showed that:
\[
R(N_0)=\frac{1}{2}\:\mathfrak{S}(N_0)N_0^2+O_A\left(\frac{N_0^2}{\log^A N_0}\right)\:, \tag{1.2}
\]
for arbitrary $A>0$, where  $\mathfrak{S}(N_0)$ is the singular series
\[
\mathfrak{S}(N_0):=\prod_{p\mid N_0}\left(1-\frac{1}{(p-1)^2}\right)\prod_{p\nmid N_0}\left(1+\frac{1}{(p-1)^3}\right)\:. \tag{1.3}
\]
The relation (1.2) implies that each sufficiently large odd integer is the sum of three primes.\\
Helfgott \cite{helfgott} recently showed that this is true for all odd $N\geq 7$. \\
Modifications of the ternary Goldbach problem are obtained by requesting that the solutions of $p_1+p_2+p_3=N$ are taken from special sets $S_i$.\\
Piatetski-Shapiro \cite{Piatetski} proved that for any fixed $c_0\in(1,12/11)$ the sequence $([n^{c_0}])_{n\in\mathbb{N}}$ contains infinitely many prime numbers. For improvements for the interval for $c_0$ cf. \cite{Deshouillers}, \cite{Heath-Brown}, \cite{Kolenski},  \cite{Rivat}. A hybrid of the theorems of Piatetski-Shapiro and Vinogradov was considered by 
Balog and Friedlander \cite{balog}. They proved that for any fixed $c_0$ with $1<c_0<21/20$ every sufficiently large odd integer $N_0$ can be  represented in the form 
$$N_0=p_1+p_2+p_3\:, \ \ \text{with $p_i=[n_i^{c_0}]$}\:,$$
for any $n_i\in\mathbb{N}$.\\
Other combinations of the sets $S_i$ were considered in \cite{balog}, \cite{jia}, \cite{Kumchev}, \cite{Teravainen}. Additionally, the authors in \cite{maier_rassias_vinogradov} proved that under the assumption of the Generalized Riemann Hypothesis each sufficiently large odd integer can be expressed as the sum of a prime and two isolated primes.\\
In the paper \cite{maier-rassias-gold}, the authors established a hybrid theorem with another type of primes: primes with missing digits. 
Numbers with restricted digits have been investigated in many papers (cf. \cite{Banks1}, \cite{Banks2}, \cite{Bourgain}, \cite{Dartyge1}, \cite{Dimitrov1},  \cite{Dimitrov2},  \cite{Dartyge2}, \cite{Drmota}, \cite{Erdos1}, \cite{Erdos2}, \cite{Konyagin}, \cite{Mauduit}). The climax of this work was the paper of Maynard \cite{maynard}, who proved the existence of infinitely many primes with restricted digits. In \cite{maynard} Maynard proved the following:\\
\textit{Let $a_0\in\{0,1,\ldots,9 \}$. Then there are infinitely many primes, whose decimal expansion does not contain the digit $a_0$.}\\
In \cite{maier-rassias-gold} we merged methods of J. Maynard \cite{maynard}, results of A. Balog and J. Friedlander \cite{balog} and the Hardy-Littlewood circle method in two variables. We proved the following:

\begin{theorem}\label{thm11}
Assume the Generalized Riemann Hypothesis (GRH). Let 
$$\gamma^*:=\frac{8}{9}+\frac{2}{3}\:\frac{\log(10/9)}{\log 10}\:(\approx 0.919\ldots)\:.$$
Let $\gamma^* \leq 1$, $c_0=1/\gamma_0$ be fixed. Let also $a_0\in\{0,1,\ldots, 9\}$.\\
Then each sufficiently large odd integer $N_0$ can be represented in the form 
$$N_0=p_1+p_2+p_3\:,$$
where the $p_i$ are of the form $p_i=[n_i^{c_0}]$, $n_i\in\mathbb{N}$, for $i=1,2$ and the decimal expansion of $p_3$ does not contain the digit $a_0$.
\end{theorem}

We now modify Theorem \ref{thm11} by replacing $p_3$, one of the two Piatetski-Shapiro primes by a prime $p$ of the form 
$$p=x^2+y^2+1\ (x, y\in \mathbb{Z})\:,$$
incorporating ideas fist used by Hooley \cite{hooley}.
Our new result is the following:

\begin{theorem}\label{thm12}
Assume the GRH. Let $\gamma^*, \gamma_0, c_0, a_0$ as in Theorem \ref{thm11}. Then each sufficiently large odd integer $N_0$ can be represented in the form 
$$N_0=p_1+p_2+p_3\:,$$
where the $p_i$ are of the form $p_2=[n_2^{c_0}]$, 
$$p_3=x_3^2+y_3^2+1\:.$$
\end{theorem}

\section{Outline and some basic definitions}

This paper is closely related to our paper \cite{maier-rassias-gold} which in turn follows Maynard \cite{maynard}. We recall the following definition from \cite{maynard}.

\begin{definition}\label{def21}
Let $a_0\in\{0,1,\ldots,9\}$, $k\in\mathbb{N}$ and let
\begin{align*}
&\mathcal{A}:=\left\{  \sum_{0\leq i\leq k} n_i 10^i \::\: n_i\in\{0,1,\ldots,9\}\setminus \{a_0\}\right\} \:,\\
&X:=10^k\:,\ \ \mathcal{B}:=\{n\leq X,\ n\in\mathbb{N}\}\:,
\end{align*}
\begin{align*}
&\mathbb{P}\ \text{the set of prime numbers},\\
&S_{\mathcal{A}}(\theta):=\sum_{a\in\mathcal{A}} e(a\theta)\:,\ \ 
S_{\mathbb{P}}(\theta):=\sum_{p\leq x}e(p\theta)\:,\ \ 
S_{\mathcal{A}\cap\mathbb{P}}(\theta):=\sum_{p\in \mathcal{A}\cap \mathbb{P}}e(p\theta)\:.
\end{align*}
Let $\mathcal{C}$ be a set of integers. We define the characteristic function $1_{\mathcal{C}}$ by
\begin{eqnarray}
1_{\mathcal{C}}(n):=\left\{ 
  \begin{array}{l l}
    1\:, & \quad \text{if $n\in \mathcal{C}$}\vspace{2mm}\\ 
    0\:, & \quad \text{if $n\not\in \mathcal{C}$}\:.\\
  \end{array} \right.
\nonumber
\end{eqnarray}
For $d\in\mathbb{N}$ we set
$$\mathcal{C}_d:=\{c\::\: cd\in\mathcal{C}\}.$$
The sifted set $\mathcal{U}(\mathcal{C}, z)$ is defined by
$$\mathcal{U}(\mathcal{C}, z):=\{c\in\mathcal{C}\::\: p\mid c\ \Rightarrow\ p>z\}\:.$$
The sieving function $S(\mathcal{C}, z)$ - the counting function of $\mathcal{U}(\mathcal{C}, z)$ - is given by
$$S(\mathcal{C}, z):=\#\mathcal{U}(\mathcal{C}, z)=\#\{c\in \mathcal{C}\::\: p\mid c\ \Rightarrow \ p>z\}\:.$$
We let
\begin{eqnarray}
w_n:=1_{\mathcal{A}}(n)-\frac{\kappa_{\mathcal{A}}\#\mathcal{A}}{\#\mathcal{B}}\:,\ \ \ 
\kappa_{\mathcal{A}}:=\left\{ 
  \begin{array}{l l}
    \dfrac{10(\Phi(10)-1)}{9\Phi(10)}\:, & \quad \text{if $(10,a_0)=1$}\vspace{2mm}\\ 
    \dfrac{10}{9}\:, & \quad \text{otherwise}\:,\\
  \end{array} \right.
\nonumber
\end{eqnarray}
$$S_d(z):=\sum_{\substack{n<X/d\\ p\mid n\ \Rightarrow\ p>z}} w_{nd}=S(\mathcal{A}_d, z)-\frac{\kappa_{\mathcal{A}} \#\mathcal{A}}{X}\: S(\mathcal{B}_d, z)\:,$$
$$1_{\mathcal{A}}(n)\ \ \text{is called the $\mathcal{A}-$part of } w_n\:,$$
$$S(\mathcal{A}_d, z)\ \ \text{is called the $\mathcal{A}-$part of } S_d(z)\:,$$
the $\mathcal{B}$-parts are defined analogously.\\
We also define the exponential sums 
$$S(\mathcal{C}, z, \theta):=\sum_{n\in\mathcal{U}(\mathcal{C},z)}e(n\theta)\:,$$
$$S_d(z, \theta):=\sum_{\substack{n<X/d\\ p\mid n\ \Rightarrow\ p>z}}w_{nd}e(n\theta)=S(\mathcal{A}_d,z,\theta)-\frac{\kappa_{\mathcal{A}}\#\mathcal{A}}{X}S(\mathcal{B}_d, z, \theta)\:,\ \ (\theta\in\mathbb{R})\:.$$
We define $X$ by $2X\leq N_0<20X$. We then define 
$$S_{c_0}:=\frac{1}{\gamma} \sum_{\frac{N_0}{2}-\frac{X}{4}<p\leq \frac{N_0}{2}-\frac{X}{8}} (\log p) p^{1-\gamma} e(p\theta)\:.$$
\end{definition}

For the proof of our modified Theorem \ref{thm12} we also have to consider the generating exponential sum for the primes of the form 
$$p=x^2+y^2+1.$$
They will be counted by multiplicity.
\begin{definition}\label{def22}
We set 
$$S_Q(\theta):=\sum_{\substack{(x,y)\in\mathbb{Z} \\ x^2+y^2+1=p\ \text{prime} \\ \frac{N_0}{2}-\frac{X}{4}<p\leq \frac{N_0}{2}-\frac{X}{8}}}  e(p\theta)\log p. $$
\end{definition}

\begin{definition}\label{def2323}
For $n\in\mathbb{N}$, let
$$r(n):=|\{(x,y)\in\mathbb{Z}^2\::\: n=x^2+y^2\}|.$$
\end{definition}

\begin{lemma}\label{lem2424}
Let $\chi$ be the non-principal Dirichlet character $\bmod 4$. For $n\in\mathbb{N}$ we have
$$r(n)=\sum_{d\mid n} \chi(d).$$
\end{lemma}
\begin{proof}
Well-known.
\end{proof}

\begin{lemma}\label{lem2525}$ $\\
(i)
$$S_Q(\theta)=\sum_{\frac{N_0}{2}-\frac{X}{4}<p\leq \frac{N_0}{2}-\frac{X}{8}} r(p-1)e(p\theta)\log p.$$
(ii)
$$S_Q(\theta)=\sum_{d\leq X} \chi(d) \sum_{\substack{\frac{N_0}{2}-\frac{X}{4}<p\leq \frac{N_0}{2}-\frac{X}{8} \\ p\equiv 1 \bmod d}} e(p\theta)\log p.$$
\end{lemma}
\begin{proof}$ $\\
(i) This follows from Definitions \ref{def22} and \ref{def2323}.\\
(ii) This follows from (i) and Lemma \ref{lem2424}.\\
We now partition $S_Q$ into three partial sums.
\end{proof}

\begin{definition}\label{def2626}
Let $D:=X^{1/2}(\log X)^{-C_0}$, $C_0>0$ to be defined later. We set:
$$S_Q^{(1)}(\theta):=\sum_{d\leq D} \chi(d)\sum_{\substack{p\in Int(N_0)\\ p\equiv 1 \bmod d}} e(p\theta) \log p,$$
$$S_Q^{(2)}(\theta):=\sum_{D<d\leq \frac{X}{D}} \chi(d)\sum_{\substack{p\in Int(N_0)\\ p\equiv 1 \bmod d}} e(p\theta) \log p,$$
$$S_Q^{(3)}(\theta):=\sum_{\frac{X}{D}<d\leq X} \chi(d)\sum_{\substack{p\in Int(N_0)\\ p\equiv 1 \bmod d}} e(p\theta) \log p,$$
with 
$$Int(N_0):=\bigg(\frac{N_0}{2}-\frac{X}{4}, \frac{N_0}{2}-\frac{X}{8}\bigg].$$
\end{definition}

\begin{lemma}\label{lem2727}
$$S_Q(\theta)=S_Q^{(1)}(\theta)+S_Q^{(2)}(\theta)+S_Q^{(3)}(\theta).$$
\end{lemma}
\begin{proof}
This follows from Lemma \ref{lem2525} (ii) and Definition \ref{def2626}.
\end{proof}

In the papers \cite{maynard} and \cite{maier-rassias-gold}, type I and type II informations are crucial. They have their origin in Harman's Sieve. We give Theorem 3.1  of Harman \cite{harman}.\\
Suppose that for any sequences of complex numbers, $a_m, b_n$, that satisfy $|a_m|\leq 1$, $|b_n|\leq 1$ we have for some $\lambda>0$, $\alpha>0$, $\beta\leq 1/2$, $M\geq 1$ and a suitable constant $Y$ that 
\[
\sum_{\substack{mn\in\mathcal{A}\\ m\leq M}} a_m=\lambda \sum_{\substack{mn\in\mathcal{B}\\ m\leq M}}a_m+O(Y) \tag{3.3.1}
\]
and
\[
\sum_{\substack{mn\in\mathcal{A}\\ X^{\alpha}\leq m\leq X^{\alpha+\beta}}} a_mb_n=\lambda \sum_{mn\in\mathcal{B}}a_mb_n+O(Y)\:, \tag{3.3.2}
\]
Let $c_r$ be a sequence of complex numbers, such that $|c_r|\leq 1$, and if $c_r\neq 0$, then 
\[
p\mid r\ \Rightarrow\ p>x^{\epsilon}\:,\ \text{for some}\ \epsilon>0.  \tag{3.3.3}
\]
Then, if $X^{\alpha}<M$, $2R<\min(X^\alpha, M)$ and $M>X^{1-\alpha}$, if $2R>X^{\alpha+\beta}$, we have
\[
\sum_{r\sim R}c_rS(\mathcal{A}_r, X^\beta)=\lambda\sum_{r\sim R} c_r S(B_r, X^\beta)+O(Y\log^3X)\:. \tag{3.3.4}
\]
The equation (3.3.1) is known as type I information, whereas  (3.3.2) is known as type II information.\\
Theorem 3.1 maybe used to obtain information on a (complicated) set $\mathcal{A}$ from a (simple) set $\mathcal{B}$. In the papers \cite{maynard} and  \cite{maier-rassias-gold}, the sets $\mathcal{A}$ and $\mathcal{B}$ are those from Definition \ref{def21}.\\An important tool in the papers \cite{maynard} and \cite{maier-rassias-gold} are Buchstab's recursions. In Maynard \cite{maynard} the $``$counting function version$"$ of the Buchstab's recursion is applied.\\
Let $u_1< u_2$. Then 
\[
S(\mathcal{C}, u_2)=S(\mathcal{C}, u_1)-\sum_{u_1<p\leq u_2}S(\mathcal{C}_p, p)\:. \tag{2.1}
\]
In our paper we use the exponential sums and obtain the following modification of (2.1):
\begin{lemma}\label{lem23}
Let $u_1<u_2$. then 
\[
S(\mathcal{C}, u_2, \theta)=S(\mathcal{C}, u_1, \theta)-\sum_{u_1<p\leq u_2}S(\mathcal{C}_p, p, \theta)\:. \tag{2.2}
\]
Like in \cite{maier-rassias-gold}, this leads to an identity
$$S_{\mathcal{A}\cap\mathbb{P}}(\theta)=\sum_{j} E_{P(\mathcal{R}_j)}(\theta)\:,$$
where the exponential sums are extended over sets $P(\mathcal{R}_j)$ linked to various polytopes $\mathcal{R}_j$, defined by linear inequalities to be satisfied by the vector
$$\left(\frac{\log p_1}{\log X}\:,\ldots, \frac{\log p_l}{\log X}  \right)\:,$$
$p_i$ being the prime factors of the numbers $n\in P(\mathcal{R}_j)$.
\end{lemma}


Like in \cite{maier-rassias-gold} we introduce another modification in our paper. Instead of considering all the integers in $\mathcal{A}$ as possible candidates for our representation of $N_0$ we now only choose the integers from a subset $\mathcal{A}^*$ of $\mathcal{A}$, which are contained in a short subinterval of $\mathcal{B}$.

\begin{definition}\label{rdef23}
Let $H\in \mathbb{N}$, $H\leq k$. For 
$$n=\sum_{j=1}^k n_j 10^j,\ \ (n_j\in \{0, \ldots, 9\})$$
we write
$$n_{H,1}:= \sum_{j=k-H+1}^k n_j 10^j=: \tilde{n}_H\cdot 10^{k-H+1}$$
and
$$n_{H, 2}:= \sum_{j=0}^{k-H} n_j 10^j\:.$$
\end{definition}

\begin{lemma}\label{rlem24}
Let $n=\tilde{n}_H\cdot 10^{k-H+1}$ as in Definition \ref{rdef23}. Then
\[
n\in \mathcal{A}\ \ \text{if and only if}\ \ \tilde{n}_H\in \mathcal{A}\ \ \text{and}\ \ n_{H,2}\in\mathcal{A}\tag{2.1}
\]
There is an integer $\tilde{n}_H\in \mathcal{A}\cap [0, 10^{H-1}]$ such that for $n^*_H:= \tilde{n}_H 10^{k-H+1}$ we have the following:
\[
\left|n^*_H-5\cdot 10^{k-1}\right|\leq \frac{3}{2}\: 10^{k-2}  \tag{2.2}
\]
and for $n_{H, 2}\in \mathcal{B}^*:= [ n_H^*, n_H^*+10^{k-H} )$ the following holds:
\[
n^*_H+n_{H,2}\in \mathcal{A}\ \ \Rightarrow\ \ n_{H, 2}\in \mathcal{A}\:.  \tag{2.3}
\]
\end{lemma}
\begin{proof}
(2.1) is obvious. To show (2.2) and (2.3) we consider the following cases:\\
\textit{Case 1}: $a_0=5$, \textit{Case 2}: $a_0=4$, \textit{Case 3}: $a_0\not\in \{4, 5\}$.\\
Specifically we have:\\
\noindent \textit{Case 1}: Let $n_{k-i}\in \{0, \ldots, 9\}\setminus\{a_0\}$ for $2\leq i\leq H-1$. Then we may take
$$n^*=4,9\cdot 10^k+\sum_{j=k-H+1}^{k-2} n_j10^j\:.$$

\noindent \textit{Case 2}: Let  $n_{k-i}\in \{0, \ldots, 9\}\setminus\{a_0\}$ for $3\leq i\leq H-1$. Then we may take
$$n^*=5,09\cdot 10^k+\sum_{j=k-H+1}^{k-3} n_j10^j\:.$$
\noindent \textit{Case 3}: The choices for $n^*$ in cases 2 and 3 are both possible.
\end{proof}


For the proof of Theorem \ref{thm12} we use the discrete ($a$-variable) circle method, which we also had used for the proof of Theorem \ref{thm11}, where we have formed the convolution $\mathcal{J}(S_{\mathcal{A}\cap\mathbb{P}}(\theta))$.\\
For an exponential sum 
$$E(\theta)=\sum v(n)e(n\theta)$$
the convolution $\mathcal{J}(E(\theta))$ was defined by
$$\mathcal{J}(E(\theta))=\frac{1}{X}\sum_{1\leq a\leq X} E\left(\frac{a}{X}\right) S_{c_0}^2\left(\frac{a}{X}\right) e\left(-N_0\frac{a}{X}\right)$$

For the proof of our new Theorem \ref{thm12} we have to replace $S_{c_0}^2$ by $S_{c_0}S_Q$ and also to consider the partition 
$$S_Q(\theta)=S_Q^{(1)}(\theta)+S_Q^{(2)}(\theta)+S_Q^{(3)}(\theta).$$

\begin{definition}\label{def24}
Let $S\subseteq [1,X]$ be a set of positive integers and $v(n)$ a sequence of real numbers. For the exponential sum
$$E(\theta):=\sum_{n\in S} v(n)e(n\theta)$$
we define for $i=1, 2, 3$: 
$$\mathcal{J}^{(i)}(E):=\frac{1}{X}\sum_{1\leq a\leq X} E\left(\frac{a}{X}\right) S_{c_0}\left(\frac{a}{X}\right) S_Q^{(i)}\left(\frac{a}{X}\right) e\left(-N_0\frac{a}{X}\right),$$
$$\mathcal{J}^{(i)}(E, \tau):=\frac{1}{X}\sum_{\frac{a}{X}\in\tau} E\left(\frac{a}{X}\right) S_{c_0}\left(\frac{a}{X}\right) S_Q^{(i)}\left(\frac{a}{X}\right) e\left(-N_0\frac{a}{X}\right),$$
for a set $\tau\subseteq [0,1]$.\\
We define 
$$\mathcal{J}(E):=\mathcal{J}^{(1)}(E)+\mathcal{J}^{(3)}(E),$$
$$\mathcal{J}(E, \tau):=\mathcal{J}^{(1)}(E, \tau)+\mathcal{J}^{(3)}(E, \tau).$$
The mean-value $M(E)$ is defined by
$$M(E):=\sum_{\substack{(m,p_2,p_3) \\ m\in S,\ p_2\in\mathbb{P}_{c_0}\cap Int(N_0),\ p_3\in S_Q\cap Int(N_0) \\ m+p_2+p_3=N_0 }} p_2^{1-\gamma}r(p_3-1)\log p_2 \log p_3 v(m)$$
\end{definition}

\begin{lemma}\label{lem25}
We have
$$\mathcal{J}(E)+\mathcal{J}^{(2)}(E):=M(E)\:.$$
\end{lemma}
\begin{proof}
This follows by orthogonality and Lemma \ref{lem2727}.\\
The transformation of the exponential sum $S_{\mathcal{A}\cap \mathbb{P}}(\theta)$ by Buchstab recursions as described in Lemma \ref{lem23} now exactly follows the transformation from \cite{maier-rassias-gold}, which in turn is based on the $``$counting function version$"$ of Maynard \cite{maynard}.\\
As in \cite{maier-rassias-gold}, it is possible to use the numerical computations of Maynard \cite{maynard}. The ($a$-variable) circle method is applied for sums, for which type I and type II information is available. These sums again give \textit{negligible} contributions, i.e. do not change the asymptotics.\\
For the major arcs estimate we find approximations of the generating series $S_{c_0}$ and $S_Q$, valid in the interals $I_{c,q}$ of the major arcs. For the Piatetski-Shapiro sum $S_{c_0}$ this has already been carried out in \cite{maier-rassias-gold}.\\
Following work of Balog and Friedlander \cite{balog}
$$S_{c_0}\left(\frac{c}{q}+\xi\right)=\frac{1}{\gamma} \sum_{\substack{\frac{N_0}{2}-\frac{X}{4}<p\leq \frac{N_0}{2}-\frac{X}{8}\\ p=[n^{1/\gamma}]}} (\log p) p^{1-\gamma} e\left(p\left(\frac{c}{q}+\xi\right)\right)$$
is approximated by $\frac{\mu(q)}{\phi(q)}M(\xi)$, where
$$M(\xi):=\sum_{\frac{N_0}{2}-\frac{X}{4}<n\leq \frac{N_0}{2}-\frac{X}{8}} e(n\xi)$$
\end{proof}

We shall use this approximation also in the present paper.\\


%



We obtain asymptotic estimates for 
$$S_Q^{(1)}\left(\frac{c}{q}+\xi\right)\ \ \text{and}\ \ S_Q^{(3)}\left(\frac{c}{q}+\xi\right)$$
for $I_{c,q}$ an interval in the major arcs, whereas the contribution of $S_Q^{(2)}\left(\frac{c}{q}+\xi\right)$ is estimated by summing over all $a/X$ and then estimating the resulting expression by sieve methods.\\
For the contribution of the other subsets in the $a$-variable circle method we closely follow \cite{maier-rassias-gold}. The $b$-variable circle method is not used in this paper, however Propositions \ref{prop63} and \ref{prop64} have been proven by it in \cite{maier-rassias-gold}.

\section{Structure of the paper}

In Section 4, we carry out the modified form of the sieve decomposition of Maynard \cite{maynard} for exponential sums - instead of counting functions - as described in the outline to reduce the proof of Theorem \ref{thm12} to the proof of three Theorems: Proposition \ref{prop41} our Type I estimate, Proposition \ref{prop42} our Type II estimate and Proposition \ref{prop43} in which the $\mathcal{A}$-part is estimated trivially.\\
In Section 5, we start evaluating the convolutions from Section 4. The range of integration is broken up in major arcs and minor arcs. The major arcs are treated as described in the outline and the minor arcs are handled by Large Sieve inequalities.\\
Proposition \ref{prop41} is handled by a method from combinatorial sieve theory, replacing the M\"obius function by functions with smaller support and Fourier Analysis to fix location and residue class of the elements. For the proofs of Propositions  \ref{prop63} and \ref{prop64} we refer to \cite{maier-rassias-gold}.\\
The dependency graph between the main statements is as follows:\\

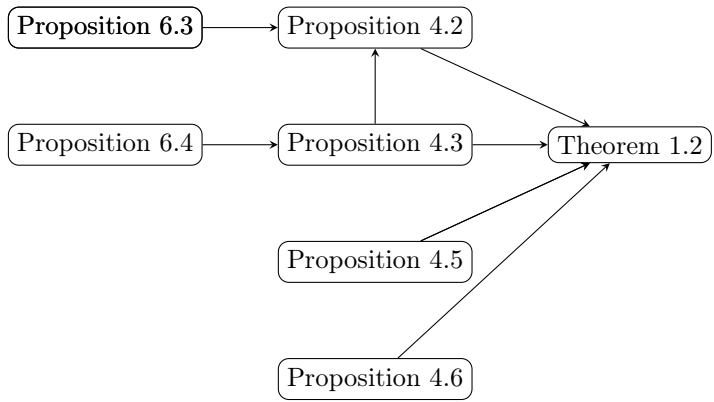
 \begin{tikzpicture}[>=stealth,every node/.style={shape=rectangle,draw,rounded corners},]
    \node (cr0) {Proposition \ref{prop63}};
        \node (cr1) [below=of cr0]{Proposition \ref{prop64}};
                 \node (c9) [right =of cr1]{Proposition \ref{prop42}};
     \node (c8) [below =of c9]{Proposition \ref{prop43}};
      \node (c10) [above =of c9]{Proposition \ref{prop41}};
            \node (cr2) [left =of c10]{Proposition \ref{prop63}};
                 \node (c12) [below =of c8]{Proposition \ref{prop4646}};

      \node (c11) [right =of c9]{Theorem \ref{thm12}};

                \draw[->] (c12) to  (c11);
                       \draw[->] (cr1) to  (c9);  
                             \draw[->] (cr2) to  (c10);  

                       \draw[->] (c9) to  (c11);
                         \draw[->] (c8) to  (c11);
                          \draw[->] (c10) to  (c11);
                            \draw[->] (c9) to  (c10);
                                 \draw[->] (c8) to  (c11);

\end{tikzpicture}

\section{Sieve decomposition and proof of Theorem \ref{thm12}}

Here we carry out the modification of Sieve Decomposition as described in the outline and reduce the proof of Theorem \ref{thm12} to the proof of Propositions \ref{prop41}, \ref{prop42},  \ref{prop43} and \ref{prop4646}.\\
Proposition \ref{prop41} deals with convolutions of exponential sums of type I, Proposition \ref{prop42} with those of exponential sums of type II, whereas Proposition \ref{prop43} gives a result in the cases in which neither type I nor type II information is available.

\begin{definition}\label{def41}
Let $\eta\in(0,1)$. Let $v(n, \eta)_{n\in\mathcal{S}}$ be a sequence of real numbers, indexed by the parameter $\eta, \mathcal{S}$ finite. The family of exponential sums
$$E(\theta; \eta):=\sum_{n\in\mathcal{S}}v(n, \eta)e(n\theta)$$
is called \textit{negligible}, if 
$$\lim_{\eta_\rightarrow 0}\lim_{k\rightarrow \infty}\frac{J(E)\log X}{|\mathcal{A}|X}=0.$$
The term $``$negligible$"$ will also be applied to an individual exponential sum $E(\theta)$ of the family $E(\theta, \eta)$.
\end{definition}

\begin{proposition}\label{prop41} (Sieve asymptotic terms)\\
Let $\epsilon>0$ be fixed. Fix an integer $l\geq 0$. Let 
$$\theta_1=\frac{9}{25}+2\epsilon\ \  \text{and}\ \  \theta_2=\frac{17}{40}-2\epsilon.$$
Let $\mathcal{L}$ be a set of $O(1)$ affine linear functions, $L\::\: \mathbb{R}^l\rightarrow \mathbb{R}$. Let
$$E_0:=E_0(\theta)=\sum^{\sim}_{X^{\theta_2-\theta_1}\leq p_1\leq \cdots \leq p_l} S^*_{p_1\cdots p_l}(X^{\theta_2-\theta_1},\theta)\:,$$
where $\sum^{\sim}$ indicates that the summation is restricted by the condition 
$$L\left( \frac{\log p_1}{\log X},\ldots, \frac{\log p_l}{\log X}\right)\geq 0\:,$$
for all $L\in\mathcal{L}$.\\
Then $E_0$ is negligible.
\end{proposition}

\begin{proposition}\label{prop42} (Type II terms)\\
Let $l\geq 1$. Let $\theta_1, \theta_2, \mathcal{C}$ be as in Proposition \ref{prop41} and let $\mathcal{I}=\{1,\ldots, l\}$ and \mbox{$j\in\{1,\ldots, l\}$,}
$$E_1(\theta):=\sum_{\substack{X^{\theta_2-\theta_1}\leq p_1\leq \cdots \leq p_l \\ X^{\theta_1}\leq \prod_{i\in I}p_i\leq X^{\theta_2}\\ p_1\cdots p_l\leq X/p_j}}^{\sim}  S^*_{p_1\cdots p_l}(p_j,\theta)\:,$$
$$E_2(\theta):=\sum_{\substack{X^{\theta_2-\theta_1}\leq p_1\leq \cdots \leq p_l \\ X^{1-\theta_2}\leq \prod_{i\in I}p_i\leq X^{1-\theta_1}\\ p_1\cdots p_l\leq X/p_j}}^{\sim}  S^*_{p_1\cdots p_l}(p_j,\theta)\:,$$
where $\sum^{\sim}$ indicates the same restriction as in Proposition \ref{prop41}.\\
Then $E_1$ and $E_2$ are negligible.
\end{proposition}

\begin{definition}\label{def42}
The Buchstab function $\omega$ is defined by the delay-differential equation
$$\omega(u)=\frac{1}{u}\:,\ 1\leq u\leq 2\:,$$
$$\omega'(u)=\omega(u-1)-\omega(u)\:,\ \ u>2\:.$$
For $\vec{p}=(p_1,\ldots,p_l)$, $p_i$ primes for $1\leq i \leq l$, let
$$Log(\vec{p})=\left(\frac{\log p_1}{\log X},\ldots, \frac{\log p_l}{\log X}\right)\:.$$
Let $\mathcal{C}$ be a set of $O(1)$ affine linear functions. Let the polytope $\mathcal{R}$ be defined by
$$\mathcal{R}:=\{(u_1, \ldots, u_l)\in[0,1]^l\::\: L(u_1,\ldots, u_l)\geq 0\ \ \text{for all}\ L\in\mathcal{C}\}\:.$$
Let 
$$P(\mathcal{R}, l):=\{\vec{p}=(p_1, \ldots, p_l)\::\: Log(\vec{p})\in \mathcal{R}\}\:$$
and
$$\Pi(\vec{p}):=p_1\cdots p_l\:.$$
We set
$$\text{Int}(N_0):=\left[ \frac{N_0}{8}-\frac{X}{8},\ \frac{N_0}{4}-\frac{X}{4} \right]\:.$$

For $y\in\mathcal{P}$ we define
$$Vol(N_0, y):=\int_{\substack{w\in \text{Int}(N_0)\\ N_0-w-y\in \text{Int}(N_0) }} 1\: dw\:.$$
$$\Gamma:=\Gamma(N_0, \mathcal{R}):=\int_{y\in{P}(R,l)}Vol(N_0,y)\: dy.$$
We define
$$\mathfrak{S}(N_0):=\prod_{p\mid N_0}\left(1-\frac{1}{(p-1)^2}\right)\prod_{p\nmid N_0}\left(1+\frac{1}{(p-1)^3}\right)\:.$$
$$\mathfrak{S}^*(N_0)=\pi \mathfrak{S}(N_0) \prod_{p\nmid N_0(N_0-1)} \left(1+\chi(p)\: \frac{p-3}{p(p^2-3p+3)}\right)\ \prod_{p\mid N_0} \left(1+\chi(p)\frac{1}{p(p-1)}\right)$$
$$\times \prod_{p\mid N_0-1}\left(1+\chi(p)\:\frac{2p-3}{p(p^2-3p+3)}\right).$$

\end{definition}

\begin{proposition}\label{prop43}
Let
$$z\::\: [0,1]^l\rightarrow \mathbb{R},\ \vec{u}=(u_1,\ldots, u_l)\rightarrow z(\vec{u})=z(u_1, \ldots, u_l)$$
be continuous on $[0,1]^l$.\\
Let 
$$E(\theta):=\sum_{\vec{p}\in P(\mathcal{R}, l)} S(\mathcal{B}_{\prod(\vec{p})}, X^{z(Log(\vec{p}))}, \theta))\:.$$
Then 
$$J(E(\theta))=\Gamma\int\cdots\int_\mathcal{R} \frac{\omega(1-u_1-\cdots - u_l)}{u_1\cdots u_l z(u_1, \ldots, u_l)}\: du_1\ldots du_l\: (1+o(1))\:.$$
\end{proposition}

\begin{proposition}\label{prop4646}
$$J^{(2)}(S_{\mathcal{A}^*\cap \mathbb{P}}(\theta))=O\left(|\mathcal{A}^*|X\:\frac{\log\log X}{\log^2 X}\right).$$
\end{proposition}

The reduction of Theorem \ref{thm12} to the proofs of Propositions \ref{prop41}, \ref{prop42} and \ref{prop43} is very similar to the analogous reduction in \cite{maier-rassias-gold}. Therefore we only describe the procedure and apply the same series of Buchstab recursions as in \cite{maier-rassias-gold}. The same sums as in \cite{maier-rassias-gold} turn out to have non-negative convolutions or are negligible - now in the sense of Definition \ref{def41} of this paper.\\
We finally obtain a result analogous to (4.14) in \cite{maier-rassias-gold}:\\
$$\mathcal{J}(S_{\mathcal{A}^*\cap \mathbb{P}}(\theta))\geq \Gamma \mathfrak{S}^*(N_0)(1+o(1))\: \frac{\kappa_{\mathcal{A}}\#\mathcal{A}^*}{\log X} 
\times (1-I_1-I_2-\cdots -I_9),$$
with $I_1+\cdots+I_9<0.996$.\\
For the definition of the integrals see \cite{maier-rassias-gold}.

\section{Sieve asymptotics}

In this section we reduce the propositions of Section 4 to other facts, some of which will be proven in this section, others will be recalled from the paper \cite{maier-rassias-gold}.\\
In the proofs of Propositions \ref{prop42} and \ref{prop43}, which deal with sums over polytopes, these polytopes will be decomposed into small boxes. We later shall recall the relevant results and definitions from \cite{maier-rassias-gold}.\\
Like in \cite{maier-rassias-gold}, also the circle method plays a crucial role. In the proofs of Propositions  \ref{prop41} and \ref{prop42} we shall apply the $a$-variable discrete circle method, whereas for Proposition \ref{prop43} we switch back to the continuous circle method. 

In both cases we split the range of summation (resp. integration) into major arcs and minor arcs.\\
As explained in the outline, we shall partition the sums $S_Q(\theta)$ into three partial sums 
$$S_Q(\theta)=S_Q^{(1)}(\theta)+S_Q^{(2)}(\theta)+S_Q^{(3)}(\theta).$$
We shall derive asymptotic approximations for $S_Q^{(1)}(\theta)$ and $S_Q^{(3)}(\theta)$ for $\theta$ in the major arcs and use them in the summation (resp. integration) over the major arcs, whereas the contribution of $S_Q^{(2)}(\theta)$ will be estimated via Parseval's equations i.e. orthogonality and sieve methods.\\
We shall also use results from \cite{maier-rassias-gold} on the asymptotics of the Piatetski-Shapiro sum $S_{c_0}(\theta)$.

\begin{lemma}\label{lem51}
Assume  the GRH. Let $c, q$ be positive integers, $(c,q)=1$. Then we have for all $\epsilon>0$:
$$\sum_{p\leq N} e\left(\left(\frac{c}{q}+\xi\right)p\right)\log p=\frac{\mu(q)}{\phi(q)} \sum_{n\leq N} e(n\xi)+O_\epsilon(N^{\frac{3}{2}+\epsilon} q(|\xi|)).\ \ (N\rightarrow \infty).$$
\end{lemma}

\begin{proof}
This is Lemma 6.1 of \cite{maier-rassias-gold}.
\end{proof}

\begin{lemma}\label{lem52}
Let $\gamma, \delta$ satisfy $0<\gamma\leq 1$, $0<\delta$ and 
$$9(1-\gamma)+12\delta<1.$$
Then, uniformly in $\alpha$, we have
$$\frac{1}{\gamma}\sum_{\substack{p<N \\ p=[n^{\frac{1}{\gamma}}]  }} e(\alpha p) p^{1-\gamma}\log p=\sum_{p<N}e(\alpha p)\log p+O(N^{1-\delta})\:,$$
where the implied constant may depend on $\gamma$ and $\delta$ only.
\end{lemma}
\begin{proof}
This is Lemma 6.2 of \cite{maier-rassias-gold}.
\end{proof}

\begin{lemma}\label{lem53}
Suppose that $\theta\in\mathbb{R}$ such that there are integers $c,q$ with $(c,q)=1$ and 
$$\left|\theta-\frac{c}{q}\right|< q^{-2}.$$
Then, for all $N\geq 2$ we have
$$\sum_{n\leq N} \Lambda(n)e(\theta n)\ll (\log N)^{\frac{7}{2}}\left(\frac{N}{q^{\frac{1}{2}}}+N^{\frac{4}{5}}+(Nq)^{\frac{1}{2}}\right).$$
\end{lemma}

\begin{proof}
This is Lemma 5.3 of \cite{maier-rassias-gold}.
\end{proof}


%

 In addition to the asymptotics of $S_{c_0}$ we also need asymptotics of $S_Q$.
 
 
\begin{lemma}\label{lem57}
Let $A>0$ be arbitrarily large, $q\leq (\log X)^{2A}$. Then we have
$$S_Q^{(1)}\left(\frac{c}{q}+\xi\right)=\sum_{d\leq D} \chi(d) c_d(c,q,1)\:\frac{1}{\phi([q,d])}\: \sum_{m_3\in \text{Int}(N_0)} e(m_3\xi)\left(1+O\left((\log X)^{-A}\right)\right)\:,$$
where
$$c_d(c,q,l)=\sum_{\substack{1\leq s\leq q \\ (s,q)=1 \\ s\equiv l \bmod ((q,d))  }} e\left(\frac{cs}{q}\right)\:.$$
\end{lemma} 
\begin{proof}
This follows from Definition \ref{def2626} by applying the GRH on the sum
$$\sum_{\substack{p\in \text{Int}(N_0) \\ p\equiv 1\bmod d }} e(p\theta)\:.$$
\end{proof}

\begin{lemma}\label{lem58}
Let $B>0$ be given. Then $C_0$ in Definition \ref{def2626} can be chosen, such that 
$$S_Q^{(3)}\left(\frac{c}{q}+\xi\right)=O\left(X(\log X)^{-B}\right)\:.$$
\end{lemma}
\begin{proof}
Here we follow \cite{maier-rassias-gold}.\\
We substitute $m_3=\frac{p_3-1}{d}$ and obtain:
$$S_Q^{(3)}\left(\frac{c}{q}+\xi\right)=\sum_{j=\pm 1} \chi(j)\: \sum_{(s,q)=1} e\left(\frac{cs}{q}\right)\ \sum_{\substack{m<0 \\ 2\mid m}}\ \sum_{\substack{p_3 \equiv s \bmod q \\ p_3 \equiv 1+jm(\bmod 4m)  }}  e(p_3\xi) \log p_3\:.$$
By partitioning the range in the inner sum into small subintervals and applying the GRH, we see that the terms $j=\pm 1$ cancel, apart from the error term claimed.\\
\end{proof}

We now give several estimates of exponential sums.
\begin{lemma}\label{lem59*}
Given $\theta\in(0,1), N\in\mathbb{N}$, then there is $q$ with $1\leq q\leq N$, such that 
$$\left| \theta-\frac{c}{q} \right|\leq \frac{1}{qN}.$$
\end{lemma}
\begin{proof}
This is the well-known approximation theorem of Dirichlet.
\end{proof}

\begin{definition}\label{def510*}
For $1\leq q\leq [X^{4/5}]+1$, $L\in \{0,\ldots, [X^{1/5}]+1\}$ we define
$$I_{c,q}(L):=\left[\frac{c}{q}-q^{-1}X^{-1}L,\ \frac{c}{q}+q^{-1}X^{-1}L  \right].$$
\end{definition}

\begin{lemma}\label{lem511*}
The intervals $I_{c,q}(L)$ cover $(0,1)$.
\end{lemma}
\begin{proof}
This follows by application of Lemma \ref{lem59*} with $N=[X^{4/5}]+1$.
\end{proof}

\begin{definition}\label{def512}
Let $\delta_0>0$ be fixed, such that
$$9(1-\gamma_0)+12\delta_0<1.$$
Then we define
$$\mathfrak{n}:=\{ \theta\in(0,1)\::\: |S_{c_0}(\theta)|\leq X^{1-\delta_0} \}$$
\end{definition}

The crucial part of the evaluation of the expression for $J^{(i)}(E)$ in Definition \ref{def24} is the contribution of the major arcs $\mathfrak{M}$. Here type I and type II information are used and one has to establish asymptotic results. Also in the case of the trivial estimate of the $\mathcal{A}$-part we need asymptotic results.\\
We now first describe how the other estimates are obtained.

\begin{lemma}\label{lem513*}
Let $E(\theta)$ be one of the exponential sums $E_j(\theta)$, $(j=0, 1, 2)$, considered in Propositions \ref{prop42}, \ref{prop43}. Then we have for $i=1, 3$:
$$\sum_{\substack{1\leq q\leq (\log X)^{C_1} \\ (c,q)=1 }} \sum_{\frac{a}{X}\in I_{c,q}\setminus I_{c,q}(L_0)} E\left(\frac{a}{X}\right) S_{c_0}\left(\frac{a}{X}\right) S_Q^{(i)}\left(\frac{a}{X}\right) e\left(-N_0\frac{a}{X}\right)$$
$$\ll |\mathcal{A}| X^2 \log X L_0^{-1}. $$
\end{lemma}
\begin{proof}
By Lemma \ref{lem52} and the GRH we have for
$$\xi\in \bigcup_{\substack{1\leq q\leq (\log X)^{C_1} \\ (c,q)=1 \\ Q<q\leq 2Q\leq (\log X)^{C_1}}} I_{c,q}\setminus I_{c,q}(L_0).$$
\[ 
S_{c_0}\left(\frac{c}{q}+\xi\right)=\frac{\mu(q)}{\phi(q)}\sum_{m\in Int} e(m\xi)+O(X^{1-\delta_0})\ll \frac{X}{QL}  \tag{5.1}
\]
From Definition \ref{def2626}, the GRH and Lemma \ref{lem58} we also obtain for $C_0$ large enough:
\[
S_Q^{(1)}\left(\frac{c}{q}+\xi\right)+S_Q^{(3)}\left(\frac{c}{q}+\xi\right) \ll \frac{X}{QL} \tag{5.2}
\]
We also have:
\[
\# \left\{ a\::\: \frac{a}{X}\in I_{c,q}(L), Q<q\leq 2Q \right\} \ll \frac{X}{QL}\tag{5.3}
\]
and 
\[
E\left(\frac{a}{X}\right) \ll |\mathcal{A}|. \tag{5.4}
\]
From (5.1), (5.2), (5.3) and (5.4) we obtain:
\[
\sum_{1\leq q\leq (\log X)^{C_1}}\sum_{\frac{a}{X}\in I_{c,q}(2L)\setminus I_{c,q}(L)} E\left(\frac{a}{X}\right) S_{c_0}\left(\frac{a}{X}\right) \left( S_Q^{(1)}\left(\frac{a}{X}\right) +  S_Q^{(3)}\left(\frac{a}{X}\right) \right) e\left(-N_0\frac{a}{X}\right)\tag{5.5}
\]
$$ \ll |\mathcal{A}| X^2 \log X L^{-1}. $$
Summation of (5.5) for $L=L_0 2^j$ gives the result of Lemma \ref{lem513*}.
\end{proof}

\begin{lemma}\label{lem514*}
$$\sum^{\sim}:=\sum_{\frac{a}{X}\in\mathfrak{n}} E\left(\frac{a}{X}\right) S_{c_0}\left(\frac{a}{X}\right)\left( S_Q^{(1)}\left(\frac{a}{X}\right) +  S_Q^{(3)}\left(\frac{a}{X}\right) \right) e\left(-N_0\frac{a}{X}\right) $$
$$\ll |\mathcal{A}|^{1/2} X^{5/2-\delta_0}$$
\end{lemma}
\begin{proof}
We apply the Cauchy-Schwarz inequality and Parseval's equation. Observing the definition of $\mathfrak{n}$ we get for $i=1,2$:
$$\sum^{\sim} \ll \left(\sum_{1\leq a\leq X} E\left(\frac{a}{X}\right)^2\right)^{1/2} X^{1-\delta_0}\left(\sum_{1\leq a \leq X} S_q^{(i)}\left(\frac{a}{X}\right)^2\right)^{1/2}$$ 
$$\ll X^{1/2} |\mathcal{A}|^{1/2} X^{1-\delta_0} X\:,$$
i.e. Lemma \ref{lem514*}.
\end{proof}

\begin{lemma}\label{lem515*}
Let $Q\leq X^{\delta_0}$, $L\leq X^{1/5}$. Then we have for $i=1, 3$:
$$\sum_{Q< q\leq 2Q}\ \sum_{c \bmod q}\ \sum_{\frac{a}{X}\in I_{c,q}(L)} E\left(\frac{a}{X}\right) S_{c_0}\left(\frac{a}{X}\right) S_Q^{(i)}\left(\frac{a}{X}\right) e\left(-N_0\frac{a}{X}\right)$$ 
$$\ll |\mathcal{A}| X^2 Q^{-1/2}.$$
\end{lemma}
\begin{proof}
We only deal with the $\mathcal{A}$-part of $E(\theta)$. We partition the intervals $I_{c,q}$ into subintervals:
$$I_{c,q}(L)=\bigcup_{j=-R_{c,q,L}}^{R_{c,q,L}}\mathcal{H}_j$$
with 
$$\mathcal{H}_j:=\left(\frac{c}{q}+q^{-1}X^{-1}2^j,\ \frac{c}{q}+q^{-1}X^{-1}2^{j+1}\right)\:,\ (j\geq 0)$$
$$R_{c,q,L}\ll \log L$$
(and analogues definition for $j<0$).\\
We write $\xi_j=q^{-1}X^{-1} 2^j$. Let 
$$z_l=\frac{a_l}{X}\in \mathcal{H}_j,\ z_l=\frac{c}{q}+\xi_l.$$
We have 
$$E(z_l)=E\left(\frac{c}{q}+\xi_l\right)+\int_{\xi_j}^{\xi_l} E'\left(\frac{c}{q}+\xi\right)\: d\xi$$
and thus 
$$|E(z_l)|\leq \left|E\left(\frac{c}{q}+\xi_j\right)\right|+\int_{q^{-1}2^j X^{-1}}^{q^{-1}2^{j+1}X^{-1}} \left|E'\left(\frac{c}{q}+\xi\right)\right|\: d\xi$$
From Lemma \ref{lem51} we obtain:
$$S_Q^{(i)}\left(\frac{c}{q}+\xi\right)\ll \frac{1}{\phi(q)} \min(X, |\xi|^{-1}),\ i=1, 2.$$
Thus we obtain
\[
\sum_{c \bmod q}\ \sum_{\frac{a}{X}\in I_{c,q}(L)} E\left(\frac{a}{X}\right) S_{c_0}\left(\frac{a}{X}\right) S_Q^{(i)} e\left(-N_0\frac{a}{X}\right) \tag{5.6}
\]
$$\ll \frac{X^2}{\phi(q)^2} \sum_{j\::\: 2^j\leq L} \min(1, q2^{-j})(1+q^{-1}2^j)$$
$$\ \ \sum_{c \bmod q}\left(  \left|E\left(\frac{c}{q}+\xi_j\right)\right|  +\int_{q^{-1}2^j X^{-1}}^{q^{-1}2^{j+1}X^{-1}} \left|E'\left(\frac{c}{q}+\xi\right)\right|\: d\xi \right)$$
Recalling the definition:
$$E(\theta):=\sum_{n\in S} v(n)e(n\theta),$$
with $v(n)=0$ for $n\not\in\mathcal{A}$ we have
\[
\sum_{c\bmod q} \left| E\left(\frac{c}{q}+\xi_j\right)\right|^2  \tag{5.7}
\]
$$=\sum_{n_1, n_2} v(n_1)v(n_2) e(n_1\xi_j) e(-n_2\xi_j)\sum_{c \bmod q} e\left( \frac{c(n_1-n_2)}{q} \right)$$
$$\ll q\sum_{s \bmod q}\ \sum_{\substack{n_1\equiv n_2 \equiv s \bmod q \\ n_1, n_2\in \mathcal{A}}} 1 \ll |\mathcal{A}|^2$$
by Lemma \ref{lem59*}.\\
We apply an analogous argument to obtain 
\[
\sum_{c \bmod q} \left| E'\left(\frac{c}{q}+\xi\right)\right|^2 \ll X^2|\mathcal{A}|^2.  \tag{5.8}
\]
From (5.6), (5.7) and (5.8) we obtain Lemma \ref{lem515*}.
\end{proof}

\begin{lemma}\label{lem516*}
Let $X^{1/5}\leq q \leq [X^{4/5}]+1$. Then the interval $I_{c,q}(L)$ is a subset of $\mathfrak{n}$.
\end{lemma}
\begin{proof}
This follows from Lemma \ref{lem53}
\end{proof}

\begin{lemma}\label{lem514}(Large sieve estimates)\\
Let
$$\mathcal{A}_1:=\left\{ \sum_{0\leq i\leq k} n_i 10^i\::\: n_i\in \{0,\ldots, 9\}\setminus \{a_0\},\ k\geq 0  \right\}$$
For $Y$ an integral power of 10, let
$$F_Y(\theta):=Y^{-\log 9/\log 10}\left| \sum_{n<Y} 1_{\mathcal{A}_1}(n) e(n\theta) \right|\:. $$
Then we have:
$$\sup_{\beta\in\mathbb{R}}\sum_{c\leq q} \sup_{|\eta|<\delta} F_Y\left(\frac{c}{q}+\beta+\eta\right)$$
$$\ll \left(1+\frac{\delta Q^2}{d}\right)\left(\left(\frac{Q^2}{d}\right)^{27/77}+\frac{Q^2}{dY^{50/77}}\right)\:.$$
\end{lemma}


\section{Local versions of Maynard's results}

In this section we reduce the propositions of Section 4 to other facts, which will be proven in later sections. The $a$-variable minor arcs estimates have been carried out in Section 5. Here we are faced with the $a$-variable major arcs part. \\
Whereas we may represent the sifted sets appearing in Propositions \ref{prop42} and  \ref{prop43} as a union of simpler sets, the set considered in Proposition \ref{prop41} is obtained by an idea related to the inclusion-exclusion principle. Its analogue in number theory in its simplest form is the Sieve of Eratosthenes-Legendre containing the M\"obius $\mu$-function $\mu(n)$:\\
Let $\mathcal{C}$ be a set of integers and $\mathcal{P}$ a set of primes then 
\begin{align*}
S(\mathcal{C}, \mathcal{P}, z)&:=\#\{n\in \mathcal{C}\::\: p\mid n,\ p\in \mathcal{P}\ \Rightarrow \ p>z\}\\
&=\sum_{n\in \mathcal{C}}\sum_{\substack{t\mid n \\ t\mid P(z)}}\mu (t)\:,\ \ \text{with}\ P(z):=\prod_{\substack{p\leq z \\ p\in \mathcal{P}}} p\:.
\end{align*}

In the theory of combinatorial sieves the M\"obius function is replaced by a function $\lambda$, having smaller support. Also in this paper we proceed in this way. The basis is the following result from Combinatorial Sieve Theory.

\begin{lemma}\label{rlem61}
Let $\kappa>0$ and $y>1$. There exist two sets of real numbers 
$$\Lambda^+=(\lambda^+_d)\ \ \text{and}\ \ \Lambda^-=(\lambda^-_d)$$
depending only on $\kappa$ and $y$ with the following properties:
\[
\lambda_1^{\pm}=1 \tag{6.1}
\]
\[
|\lambda_d^{\pm}| \leq 1\:,\ \ \text{if}\ 1\leq d< y \tag{6.2}
\]
\[
\lambda_d^{\pm}=0\:,\ \ \text{if}\ d\geq  y 
\]
and for any integer $n>1$,
\[ 
\sum_{d\mid n} \lambda_d^-  \leq 0  \leq \sum_{d\mid n} \lambda_d^+\:.  \tag{6.3}
\]
Moreover, for any multiplicative function $g(d)$ with $0\leq g(p)<1$ and satisfying the dimension conditions 
$$\prod_{w\leq p< z} (1-g(p))^{-1}\leq \left(\frac{\log z}{\log w}\right)^\kappa \left(1+\frac{\kappa}{\log w}\right)$$
for all $2\leq w< y$, we have
$$\sum_{d \mid P(z)} \lambda_d^\pm g(d)=\left(1+O\left(e^{-s}\left(1+\frac{\kappa}{\log z}\right)^{10}\right)\right)\prod_{p< z} (1-g(p))\:,$$
where $P(z)$ denotes the product of all primes $p<z$ and $s=\log y/\log z$. The implied constants depend only on $\kappa$.
\end{lemma}
\begin{proof}
This is the Fundamental Lemma 6.3 of \cite{IKW}.
\end{proof}

A special role in the consideration of the sifted set in Proposition \ref{prop41} is played by the prime factors of 10, $p=2$ and $5$. To simplify things we consider the subsets
$$\mathcal{U}^{*'}:=\{a\in \mathcal{U}^*\::\: (a, 10)=1\}$$
and
$$\mathcal{B}^{*'}:=\{ b\in \mathcal{B}^*\::\: (b, 10)=1\}$$

\begin{definition}\label{def62}
Let $\lambda$ be an arithmetic function, $z\geq 1$. For a set of $\mathcal{C}$ of integers we define 
$$S(\mathcal{C}, z, \theta, \lambda):=\sum_{n\in \mathcal{C}} e(n\theta)\left(\sum_{\substack{t\mid n \\ t\mid P(z)}} \lambda(t)\right)$$
\end{definition}
We are now ready for the statement of 

\begin{proposition}\label{prop63}
Let $\epsilon>0$, $0<\eta_0\leq \theta_2-\theta_1$, $l=l(\eta_0)$ be fixed. Let $\mathcal{L}$, and the summation condition $\sum^\sim$ be as in Proposition \ref{prop41}, $q\leq Q_0$, $(c, q)=1$. Let $\lambda^\pm$ satisfy the properties of Lemma \ref{rlem61} with $y=X^{(\eta_0^{1/2})}$ and let $\lambda^\pm(t)=0$, if $(t, 10)>1$. Then we have for $\lambda=\lambda^-$ or $\lambda^+$:
$$\sum_{X^{\eta_0}\leq p_1\leq \cdots \leq p_l}^\sim \left( S\left(\mathcal{A}^*_{p_1 \ldots p_l}, X^{\eta_0}, \frac{c}{q}, \lambda\right)-
\frac{\kappa_{\mathcal{A}}\#\mathcal{A}^*}{\#\mathcal{B}^*}  S\left(\mathcal{B}^*_{p_1 \ldots p_l}, X^{\eta_0}, \frac{c}{q}, \lambda\right)  \right)$$
$$=O\left((\#\mathcal{A}^*)(\log X)^{-A}\right)\:.$$
\end{proposition}

\begin{proposition}\label{prop64}(Type II terms, local version)\\
Let $\epsilon, \eta_0, l, \mathcal{L}, \sum^\sim, q, c, t$ be as in Proposition \ref{prop63}. Then we have
$$\sum^\sim_{\substack{X^{\eta_0}\leq p_1\leq \cdots \leq p_l \\ X^{\theta_1}\leq \prod_{i\in \mathcal{I}} p_i\leq X^{\theta_2}\\ p_1\cdots p_l\leq X/p_j}} \left( S\left(\mathcal{A}^*_{p_1 \ldots p_l}, p_j, \frac{c}{q}\right)-\frac{\kappa_{\mathcal{A}}\#\mathcal{A}^*}{\#\mathcal{B}^*}  S\left(\mathcal{B}^*_{p_1 \ldots p_l}, p_j, \frac{c}{q}\right)  \right) $$
$$=O\left((\#\mathcal{A}^*)(\log X)^{-A}\right)\:.$$
and
$$\sum^\sim_{\substack{X^{\eta_0}\leq p_1\leq \cdots \leq p_l \\ X^{1-\theta_2}\leq \prod_{i\in \mathcal{I}} p_i\leq X^{1-\theta_1}\\ p_1\cdots p_l\leq X/p_j}} \left( S\left(\mathcal{A}^*_{p_1 \ldots p_l}, p_j, \frac{c}{q}\right)-\frac{\kappa_{\mathcal{A}}\#\mathcal{A}^*}{\#\mathcal{B}^*}  S\left(\mathcal{B}^*_{p_1 \ldots p_l}, p_j, \frac{c}{q}\right)  \right) $$
$$=O\left((\#\mathcal{A}^*)(\log X)^{-A}\right)\:.$$
\end{proposition}


We now complete the reduction of Proposition \ref{prop41} to Proposition \ref{prop63} and Proposition \ref{prop42}. Here we closely follow \cite{maier-rassias-gold}. The main change consists  in replacing the factor $S_{c_0}^2$ by $S_{c_0}S_{Q}^{(i)}$ $(i=1, 3)$. We first prove a modification containing the weights $\lambda^\pm$.

\begin{lemma}\label{rlem616}
Let $\lambda^\pm$ satisfy the properties of Lemma \ref{rlem61} and $\lambda^\pm(t)=0$, if \mbox{$(t, 10)>1$.} Let
$$E_{0, \mathcal{A}^*, \lambda}(\theta)=\sum_{x^{\eta_0}\leq p_1\leq \cdots \leq p_l}^\sim S(\mathcal{A}_{p_1\cdots p_l}^*, X^{\eta_0}, \theta, \lambda)$$
$$E_{0, \mathcal{B}^*, \lambda}(\theta)=\sum_{x^{\eta_0}\leq p_1\leq \cdots \leq p_l}^\sim S(\mathcal{B}_{p_1\cdots p_l}^*, X^{\eta_0}, \theta, \lambda)$$
Then for $\lambda=\lambda^-$ or $\lambda^+$, $i=1$ or $3$, we have:
$$\frac{1}{X}\sum_{1\leq a\leq X} \left( E_{0, \mathcal{A}^*, \lambda}\left(\frac{a}{X}\right)-\kappa_{\mathcal{A}}\frac{\#\mathcal{A}^*}{\#\mathcal{B}^*} E_{0, \mathcal{B}^*, \lambda}\left(\frac{a}{X}\right)\right)$$
$$S_{c_0}\left(\frac{a}{X}\right)S_Q^{(i)}\left(\frac{a}{X}\right)e\left(-N_0 \frac{a}{X}\right)=O\left(\# \mathcal{A}^* X(\log X)^{-A}  \right).$$
\end{lemma}

\noindent\textit{Proof of Lemma \ref{rlem616} assuming Proposition \ref{prop63} }

The first step in the approximation of $E_{0, \mathcal{A}^*, \lambda}(\theta)$ inside $I_{c,q}$ consists in the replacement of the variable factors $e\left( n\left(\frac{c}{q}+\xi\right)\right)$ by $e\left(n_0\frac{c}{q}\right)e(n\xi)$, where $n_0$ is the midpoint of the interval $\mathcal{B}^*$. We set
$$\tilde{\lambda}(n)=\sum_{t\mid n} \lambda (t)$$
and obtain:
\begin{align*}
E_{0, \mathcal{A}^*, \lambda}\left(\frac{c}{q}+\xi\right)& =\sum_{X^{\eta_0}\leq p_1\leq \cdots\leq p_l}^\sim S\left(\mathcal{A}^*_{p_1\ldots p_l}, X^{\eta_0}, \left(\frac{c}{q}+\xi\right), \lambda\right)   \tag{6.4}\\
&=\sum_{X^{\eta_0}\leq p_1\leq \cdots\leq p_l}^\sim e(n_0\xi) \sum_{n\in\mathcal{B}^*}  1_{\mathcal{A}^*_{p_1\ldots p_l}}(n) e\left(n\frac{c}{q}\right) \tilde{\lambda}(n)\\
&\ \ +\sum_{X^{\eta_0}\leq p_1\leq \cdots\leq p_l}^\sim e(n_0\xi) \sum_{n\in\mathcal{B}^*}  1_{\mathcal{A}^*_{p_1\ldots p_l}}(n) e\left(n\frac{c}{q}\right) (e(n\xi)-e(n_0\xi))\tilde{\lambda}(n)\\
&=:  E_{0, \mathcal{A}^*, \lambda}^{(1)}+ E_{0, \mathcal{A}^*, \lambda}^{(2)}
\end{align*}

We have

\[
E_{0, \mathcal{A}^*, \lambda}^{(1)} = e(n_0\xi)\sum_{X^{\eta_0}\leq p_1\leq \cdots \leq p_l}^\sim  S\left(\mathcal{A}^*_{p_1\ldots p_l}, X^{\eta_0}, \frac{c}{q}, \lambda\right)    \tag{6.5}
\]
and
\[
E_{0, \mathcal{A}^*, \lambda}^{(2)} = O\left(\sum_{X^{\eta_0}\leq p_1\leq \cdots \leq p_l}^\sim \ \sum_{n\in\mathcal{B}^*} 1_{\mathcal{A}^*_{p_1\ldots p_l}}(n)\: |n-n_0|\:|\xi| \:\tilde{\lambda}(n)\right) \:.     \tag{6.6}
\]

We obtain an analogous decomposition for the $\mathcal{B}^*$-part:
$$E_{0, \mathcal{B}^*, \lambda}:=E^{(1)}_{0,\mathcal{B}^*, \lambda }+E_{0, \mathcal{B}^*, \lambda}^{(2)} \:.$$

From (6.4), (6.5), and (6.6) we obtain:
\[ 
\left| E_{0, \mathcal{A}^*, \lambda}\left(\frac{c}{q}+\lambda\right)- \kappa_{\mathcal{A}}\frac{\#\mathcal{A}^*}{\#\mathcal{B}^*} E_{0, \mathcal{B}^*, \lambda}\left(\frac{c}{q}+\lambda\right)\right|
 \tag{6.7}
\]
$$\leq \left| E_{0, \mathcal{A}^*, \lambda}\left(\frac{c}{q}\right)-\kappa_{\mathcal{A}}\frac{\#\mathcal{A}^*}{\#\mathcal{B}^*} E_{0, \mathcal{B}^*, \lambda}\left(\frac{c}{q}\right)\right|+|E_{0, \mathcal{A}, \lambda}^{(1)}|+\kappa_{\mathcal{A}}\frac{\#\mathcal{A}^*}{\#\mathcal{B}^*} |E_{0, \mathcal{A}, \lambda}^{(2)}|.$$
By Proposition \ref{prop63} we have:
\[
 \left| E_{0, \mathcal{A}^*, \lambda}\left(\frac{c}{q}\right)-\kappa_{\mathcal{A}}\frac{\#\mathcal{A}^*}{\#\mathcal{B}^*} E_{0, \mathcal{B}^*, \lambda}\left(\frac{c}{q}\right)\right|=O\left( \#\mathcal{A}(\log X)^{-A}\right)\:.  \tag{6.8}
\]
We also have:
\[
e(n\xi)-e(n_0\xi)=O(|n-n_0|\:|\xi|)\:.  \tag{6.9}
\]
Additionally, we have
\[
\sum_{X^{\eta_0}\leq p_1\leq \cdots\leq p_l}^{\sim}\: \sum_{n\in \mathcal{A}^*_{p_1\cdots p_l}} \sum_{t\mid n} \lambda(t)\ll \sum_{\substack{X^{\eta_0}\leq p_1\leq \cdots\leq p_l \\ t\leq X^\rho}}\ \sum_{n\in \mathcal{A}_{[p_1\cdots p_l, t]}} 1 \tag{6.10}
\]

From Lemma \ref{lem514} we have
\[
\# \mathcal{A}^*_{[p_1\ldots p_l, t]} = O\left( \frac{\#\mathcal{A}^*}{[p_1\cdots p_l, t]} \right)  \tag{6.11}
\]
The major arcs estimate for Lemma \ref{rlem616} can now be concluded:\\
From (6.7), (6.8), (6.9), (6.10) we obtain:
$$\frac{1}{X}\sum_{q\leq Q_0}\sum_{(c,q)>1}\sum_{\frac{a}{X}\in I_{c,q}(L_0)} \left( E_{0, \mathcal{A}^*, \lambda}\left(\frac{a}{X}\right) - \frac{\kappa_{\mathcal{A}}\#\mathcal{A}^*}{\#\mathcal{B}^*} E_{0, \mathcal{B}^*, \lambda}\left(\frac{a}{X}\right) \right) S_{c_0}\left(\frac{a}{X}\right) S_Q^{(i)}\left(\frac{a}{X}\right) e\left(-N_0\frac{a}{X}\right)$$
$$ = O\left(  (\#\mathcal{A}^*) X(\log X)^{-A}\right)\:.$$
The minor arcs estimates are now obtained by treating the $\mathcal{A}$-part $E_{0, \mathcal{A}^*, \lambda}(\theta)$ and the $\mathcal{B}$-part $E_{0, \mathcal{B}^*, \lambda}(\theta)$ separately. The estimates are easily carried out by the application of Lemmas \ref{lem511*}, \ref{lem513*}, \ref{lem514*}, \ref{lem515*} and \ref{lem516*}.

\textit{Conclusion of the proof of Proposition \ref{prop41}}

\noindent The next step consists in replacing the functions $\lambda^{\pm}$ from Lemma \ref{rlem616} by the M\"obius function $\mu$, thus obtaining the exponential sum $E_0(\theta, \eta)$ from Proposition \ref{prop41}. We set
$$\mathcal{U}^{*'}:=\{ m\in \mathcal{U}^*\::\: (m, 10)=1\}$$
$$\mathcal{B}^{*'}:=\{ n\in \mathcal{B}^*\::\: (n, 10)=1\}$$
and observe that
$$S(\mathcal{U}^{*'}, X^{\eta_0}, \theta, \mu)=S(\mathcal{U}^*, X^{\eta_0}, \theta, \mu)$$
and
$$S(\mathcal{B}^{*'}, X^{\eta_0}, \theta, \mu)=S(\mathcal{B}^*, X^{\eta_0}, \theta, \mu)\:.$$
because of the condition $\lambda(t)=0$ for $(t, 10)>1$ we have:
\[
J(E(\mathcal{A}^*, X^{\eta_0}, \lambda^-))\leq J(E(\mathcal{A}^*, X^{\eta_0}, \mu))\leq J(E(\mathcal{A}^*, X^{\eta_0}, \lambda^+)) \tag{6.12}
\]
and
\[
J(E(\mathcal{B}^*, X^{\eta_0}, \lambda^-))\leq J(E(\mathcal{B}^*, X^{\eta_0}, \mu))\leq J(E(\mathcal{B}^*, X^{\eta_0}, \lambda^+)). \tag{6.13}
\]
We now apply Lemma \ref{rlem61} with 

\begin{eqnarray}
g(p):=\left\{ 
  \begin{array}{l l}
   0\:, & \quad \text{if $p\in \{2, 5\}$}\vspace{2mm}\\ 
    1/p\:, & \quad \text{otherwise}\:,\\
  \end{array} \right.
\nonumber
\end{eqnarray}

and obtain:\\
For all $\epsilon>0$ there is an $\eta^*$, such that
$$\limsup_{k\rightarrow\infty} \frac{|J(E(\theta, \eta^*))|\log X}{|\mathcal{A}^*|X}< \epsilon\:,\ \ \text{for}\ \eta^*\leq \eta_0\:.$$

We still have to pass from $X^{\eta_0}$ to $X^{\theta_2-\theta_1}$. We modify the analysis in \cite{maynard}, p. 156:\\
Given a set $\mathcal{C}$ and an integer $d$ we let
$$T_m(\mathcal{C}; d, \theta):=\sum_{\substack{X^\eta\leq p_m'\leq \cdots \leq p_1'\leq X^\theta \\ d{p_1'\cdots p_m'}\leq X^{\theta_1}}} S(\mathcal{C}_{p_1'\cdots p_m'}, X^\eta, \theta) $$
$$U_m(\mathcal{C}; d, \theta):=\sum_{\substack{X^\eta\leq p_m'\leq \cdots \leq p_1'\leq X^\theta \\ d{p_1'\cdots p_m'}\leq X^{\theta_1}}} S(\mathcal{C}_{p_1'\cdots p_m'}, p_m'X^\eta, \theta) $$
$$V_m(\mathcal{C}; d, \theta):=\sum_{X^\eta< p_m'\leq \cdots \leq p_1'\leq X^\theta } S(\mathcal{C}_{p_1'\cdots p_m'}, p_m', \theta) .$$
Buchstab's identity shows that 
$$U_m(\mathcal{C}; d, \theta)=T_m(\mathcal{C}; d, \theta)-U_{m+1}(\mathcal{C}; d, \theta)-V_{m+1}(\mathcal{C}; d, \theta)$$
The $T_m$-terms are now handled by Lemma $\ref{rlem616}$, whereas the $V_m$-terms are reduced to Proposition \ref{prop42}.\\
Proposition \ref{prop41} now has been reduced to Proposition \ref{prop63} and Proposition \ref{prop42}.

\textit{Proof of Proposition \ref{prop42} assuming Proposition \ref{prop64}}

\noindent We restrict ourselves to $E_1(\theta, \eta_0)$, since the case of $E_2(\theta, \eta_0)$ is completely analogous. As in the proof of Proposition \ref{prop41} we replace the variable factors $e(n(\frac{c}{q}+\xi))$ by $e(n_0\xi)e(n\frac{c}{q})$ with $n_0\in\mathcal{B}^*$.\\
We obtain
$$\sum_{0, \mathcal{A}^*}e\left(\frac{c}{q}+\xi\right)=\sum_{n}1_{U(\mathcal{A}^*_{p_1\ldots p_l}, p_j)}(n)\:e\left(n\left(\frac{c}{q}+\xi\right)\right)=\Sigma_{0, \mathcal{A}^*}^{(1)}  +\Sigma_{0, \mathcal{A}^*}^{(2)}  $$
with
$$\Sigma_{0, \mathcal{A}^*}^{(1)}= e(n_0\xi) \sum_{n}1_{U(\mathcal{A}^*_{p_1\ldots p_l}, p_j)}(n) e\left(n\frac{c}{q}\right)$$
and
$$\Sigma_{0, \mathcal{A}^*}^{(2)}= e(n_0\xi) \sum_{n}1_{U(\mathcal{A}^*_{p_1\ldots p_l}, p_j)}(n) ( e(n\xi)-e(n_0\xi)).$$
An analogous decomposition holds for
$$\sum_{0, \mathcal{B}^*}e\left(\frac{c}{q}+\xi\right)\:.$$
The claim of Proposition \ref{prop42} now follows quite analogously to the proof of Lemma \ref{rlem616}. We use Proposition \ref{prop64} for the estimate of
$$\Sigma_{0, \mathcal{A}^*}^{(1)}  -\frac{\kappa_{\mathcal{A}}\#\mathcal{A}^*}{\#\mathcal{B}^*}  \Sigma_{0, \mathcal{B}^*}^{(1)}\:, $$
where for the other major arcs contribution we again use the estimate 
$$|e(n\xi)-e(n_0\xi)|=O(|n-n_0|\:|\xi|)\:.$$
The minor arcs estimates follow again by the application of Lemmas \ref{lem511*}, \ref{lem513*}, \ref{lem514*}, \ref{lem515*} and \ref{lem516*}.

\textit{Proof of Proposition \ref{prop43}}

\noindent We first deal with  the $a$-variable major arcs contribution:\\
Let $1\leq q\leq Q_0$, $(c, q)=1$, $\eta=q^{-1}X^{-1}L_0$. By Lemma \ref{rlem616} and the GRH we have for $|\xi|\leq \eta$:
$$S_{c_0}\left(\frac{c}{q}+\xi\right)=\frac{\mu(q)}{\phi(q)}\:\sum_{m\in Int(N_0)} e(m\xi)+O(X^{1-\delta_0}).$$

We now approximate $E\left(\frac{c}{q}+\xi\right)$. For $n\in \mathcal{U}\left( \mathcal{B}^*_{\prod(\vec{p})}, X^{z(Log \vec{p})}\right)$ we write
$$n=p_1\cdots p_l\cdot m\ \ \text{and}\ \ m=q_1\cdots q_r$$
with
$$X^{z(Log \vec{p})} \leq q_1< q_2<\cdots <q_v.$$
By partitioning the range of the $p_k$  and the $q_j$ into intervals and using GRH we see
\begin{align*}
\mathcal{U}(q, s)&:=\# \left\{ n\in \mathcal{U}\left( \mathcal{B}^*_{\prod(\vec{p})}, X^{z(Log \vec{p})}\right),\ n\equiv s \bmod q\right\}  \tag{6.14}\\
&=\mathcal{U}(q, s_0)\left(1+O((\log X)^{-A})\right)\:,
\end{align*}
for any $s_0$ with $(q, s_0)=1$, i.e. $\mathcal{U}(q, s)$ is asymptotically independent of $s$.\\
From (6.14) we obtain:
\begin{align*}
E\left(\frac{c}{q}+\xi\right)&=\sum_{\substack{s \bmod q \\ (s, q)=1}} e\left(\frac{sc}{q}\right) \sum_{n\in \mathcal{U}(q, s)} e(n\xi)(1+O(\log X)^{-C_4})\\
&=\frac{\mu(q)}{\phi(q)} \sum_{\vec{p}\::\: Log\: \vec{p}\:\in\: \mathcal{R}}\ \  \sum_{m\:\in\: \mathcal{U}(\mathcal{B}_{\prod(\vec{p})}, X^{z(Log\: \vec{p})})}  (1+O(\log X)^{-C_4}).
\end{align*}

We obtain

\begin{align*}
&\int_{-\eta}^{\eta} E\left(\frac{c}{q}+\xi\right) S_{c_0}\left(\frac{c}{q}+\xi\right)S_{Q}^{(i)}\left(\frac{c}{q}+\xi\right) e\left(-N_0\left(\frac{c}{q}+\xi\right)\right) d\xi  \tag{6.15}\\
&=\frac{\mu(q)^3}{\phi(q)^3} e\left(-N_0\frac{c}{q}\right) \int_{-1/2}^{1/2} E(\xi)\sum_{(n_1, n_2)\in Int} e(n_1\xi+n_2\xi) e(-N_0\xi) d\xi\\
&=\frac{\mu(q)^3}{\phi(q)^3}\ e\left(-N_0\frac{c}{q}\right) \sum_{\vec{p}\::\: Log\: \vec{p}\:\in\: \mathcal{R}} \#\Big\{(m, n_1, n_2)\::\: m\in \mathcal{U}\left( \mathcal{B}^*_{\prod(\vec{p})}, X^{z(Log\: \vec{p})}\right),\\
&\ \ \ \  n_i\in Int, m+n_1+n_2=N_0\Big\}\:.
\end{align*}

We write $m=p_1\cdots p_l\cdot h$ with $m=q_1\cdots q_v$. By the well-known connection between the Buchstab function and the number of integers free of small prime factors, we have:
\begin{align*}
&\#\left\{ h\::\: h\in \frac{\# \mathcal{B}^*}{p_1\cdots p_l}\::\: p(h\Rightarrow p\geq z(Log(\vec{p}))\right\} \tag{6.16}\\
&=\frac{\# \mathcal{B}^*}{p_1\cdots p_l} \omega\left( \frac{\log (X/p_1\cdots p_l)}{Log(\vec{p})}\right) \frac{1}{\log X} (1+o(1))\:.
\end{align*}
The function
$$M(q)=\sum_{(c, q)=1}e\left(-N_0\frac{c}{q}\right)$$
is a multiplicative function of $q$: We obtain the singular series $\mathfrak{S}(N_0)$.\\
From (6.14), (6.15), (6.16) we obtain the major arcs contribution:
\begin{align*}
&\sum_{q\leq Q_0} \sum_{(c, q)=1} \int_{\frac{c}{q}-\eta}^{\frac{c}{q}+\eta} E(\xi)S_{c_0}(\xi)S_{Q}^{(i)}(\xi) e(-N_0\xi)d\xi\\
&= \frac{X(\# \mathcal{B}^*)}{4\log X} \mathfrak{S}_0(N_0) \int\cdots\int_\mathcal{R} \frac{\omega(1-u_1-\cdots - u_l)}{u_1\cdots u_l z(u_1, \ldots, u_l)} du_1\cdots du_l (1+o(1))\:.
\end{align*}

The proof of Proposition \ref{prop43} is complete by application of Lemmas \ref{lem511*}, \ref{lem513*}, \ref{lem514*}, \ref{lem515*} and \ref{lem516*}.

%


\section{Proof of Proposition \ref{prop4646}}

By orthogonality we have, that 
\[
\frac{1}{X}\sum_{1\leq a\leq X} S_{\mathcal{A}\cap\mathbb{P}}\left(\frac{a}{X}\right) S_{c_0}\left(\frac{a}{X}\right) S_Q^{(2)}\left(\frac{a}{X}\right) e\left(-N_0\frac{a}{X}\right) \tag{7.1}
\]
$$=\sum_{\substack{(p_1, p_2, p_3) \\ p_1+p_2+p_3=N_0 \\ p_1\in \mathcal{A^*}\cap \mathbb{P}, p_2\in\mathbb{P}_{c_0}\cap Int, p_3\in \mathbb{P}_Q\cap Int}}\ \sum_{\substack{d\mid p_3-1\\ D<d\leq \frac{X}{D}}} \chi(d) $$
$$=\sum_{D<d\leq \frac{X}{D}}\ \ \ \sum_{\substack{(p_1, p_2, p_3), p_1\in \mathcal{A}\cap \mathbb{P}, p_2\in\mathbb{P}_{c_0}\cap Int, p_3\in\mathbb{P}_Q\cap Int \\ p_1+p_2+p_3=N_0 \\ p_3\equiv 1\bmod d}} 1$$
$$=:\sum_{D<d\leq \frac{X}{D}}\ \sum (d),\ \ \text{say}.$$

For the evaluation of the inner sum we partition the primes $p_i$ into residue classes $\bmod d$. We write  \mbox{$p_i=l_i+t_id$}. From the conditions
$$p_1+p_2+p_3=N_0,\ \ p_3\equiv 1 \bmod d$$
we ontain:
\[
p_3\equiv N_0-l_1-l_2 \bmod d\:,\tag{7.2}
\]
\[
l_1+l_2\equiv N_0-1 \bmod d\:.\tag{7.3}
\]

We now collect the contributions to the inner sum $\sum (d)$ in (7.1) for a fixed prime $p_2\in\mathbb{P}_{c_0}$, which by (7.2), (7.3) means for a fixed triplet $(l_2, t_2, d)$. We set
$$h_1(t_1)=h_1(t_1; l_2, t_2, d)=l_1+t_1d$$
$$h_3(t_1)=h_3(t_1; l_2, t_2, d)= 1-(t_1+t_2)d.$$
We obtain
$$\sum(d)=\sum_{(l_2, t_2)} N(l_2, t_2, d),$$
where 
$$N(l_2, t_2, d)=\#\{ t_1\::\: h_1(t_1)\in \mathcal{A^*}\cap \mathbb{P}, h_3(t_1)\in\mathbb{P} \}.$$

To obtain an upper bound for $N(l_2, t_2, d)$, we apply Lemma \ref{rlem61} with 
$$C=\{ h_1(t_1) h_3(t_1)\::\: h_1(t_1)\in \mathcal{A}^*  \}.$$
Let $\Lambda^+, \Lambda^{-}$ be chosen as in Lemma \ref{rlem61}. For the application of Lemma \ref{rlem61} we need an estimate for
$$\#\mathfrak{C}_e=|\{h_1(t_1) h_3(t_1) \equiv 0 \bmod e,\ h_1(t_1)\in \mathcal{A}^* \}|.$$
Then we obtain
$$S(\mathcal{C}, z)\leq \sum_{e\mid P(z)} \lambda^+(e) |\mathfrak{C}_e|$$
As in the proof of Proposition \ref{prop41} we approximate the characteristic function of $\mathcal{A}^*$ by Fourier series and obtain
$$\sum_{t_1: h_1(t_1)\in \mathcal{A}_e}=\frac{1}{[d, e_1]}\sum_{a\in \mathcal{A}^*}\sum_{m \bmod [d,e]} \sum_{n=-\infty}^\infty \alpha_n e\left( mn\:\frac{(a-\tilde{l})}{[d,e]} \right)\:,$$
where $\tilde{l}$ is determined by the system
$$\tilde{l}\equiv l_1 \bmod d$$
$$\tilde{l} \equiv 0 \bmod e.$$
The non-constant terms of the Fourier series give a negligible error.\\
We may thus apply Lemma \ref{rlem61} with the multiplicative function $g(d)$ defined by 
$$g(r)=\#\{ t_1 \bmod r\::\: t_1 \equiv -l_1 d^{-1} \bmod r, $$
$$ t_1\equiv d^{-1} (N_0-(l_1+l_2)-t_2d) \bmod r  \}$$
for prime $r$.\\
By summing over all triplets $(l_2, t_2, d)$ for which $p_2=t_2+l_2d_2\in \mathbb{P}_{c_0}$. This gives the result of Proposition  \ref{prop4646}.


%
%
%
%

\end{document}